\documentclass[11pt,a4paper]{amsart}

\usepackage{epsfig}
 
\def\Conv {{\rm \mbox{Conv}}} 
\def\Vol {{\rm \mbox{Vol}}} 
\def\dist {{\rm \mbox{dist}}} 
\def\CM{{\rm \mbox{CM}}}

\def\C{\mathbb{C}}

\def \N{\mathbb{N}}

\def \Q{\mathbb{Q}}
\def \R{\mathbb{R}}
\def \Z{\mathbb{Z}}

\def\ov#1{{\overline{#1}}}

\newtheorem{defn}{Definition}[section]
\newtheorem{lem}[defn]{Lemma}

\newtheorem{prop}[defn]{Proposition}
\newtheorem{thm}[defn]{Theorem}

\newtheorem{rem}[defn]{Remark}

           {\vspace{3.3mm}
           \noindent{\bf #1}\it}%
           {\vspace{3.3mm}}

\newenvironment{Proof}[1]{{\it #1}}{\hfill\mbox{$\Box$} }
      
\begin{document}

\title[The Cayley-Menger determinant is irreducible for~$n\geq3$]
{The Cayley-Menger determinant is irreducible
  for~$\boldsymbol{{\mathnormal n} \geq 3}$}

\author[Carlos~D'Andrea]{Carlos~D'Andrea}

\address{Department of Mathematics, 
University of California at Berkeley, 
Berkeley CA~94720, USA.}

\email{cdandrea@math.berkeley.edu}

\urladdr{http://math.berkeley.edu/\~{}cdandrea/}

\author[Mart{\'\i}n~Sombra]{Mart{\'\i}n~Sombra}

\address{Universit{\'e} de Lyon 1, 
Laboratoire de Math{\'e}matiques 
Appliqu{\'e}es de Lyon, 
21 avenue Claude~Bernard,
69622 Villeurbanne Cedex, 
France.}

\email{sombra@maply.univ-lyon1.fr}

\urladdr{http://maply.univ-lyon1.fr/\~{}sombra/} 

\date{\today} 



\thanks{C.~D'Andrea was supported by a Miller Research Fellowship.}

\begin{abstract}
We prove that the Cayley-Menger determinant of an $n$-dimensional simplex is an absolutely irreducible polynomial for $n\geq3.$
We also study the irreducibility of polynomials associated to 
related geometric constructions. 
\end{abstract}
\maketitle

\overfullrule=1mm

\vspace{-8mm} 


\setcounter{section}{1}

\section*{}

Let $\{d_{ij} \, : \ 0\le i < j \le n\}$ be a set of 
$\displaystyle \frac{n\, (n+1)}{2}$ variables and consider the square 
$(n+2) \times (n+2)$  matrix
\begin{equation} \label{pn} 
\CM_n:= 
 \left[ \begin{array}{cccccc}
0 & 1         & 1        & 1        & \cdots & 1        \\[2mm]
1 & 0         & d_{0 1}^2 & d_{0 2}^2 & \cdots & d_{0 n}^2 \\[2mm] 
1 & d_{0 1}^2  & 0        & d_{1 2}^2 & \cdots & d_{1 n}^2 \\[2mm] 
1 & d_{0 2}^2  & d_{1 2}^2 & 0        & \cdots & d_{2 n}^2 \\[2mm] 
\vdots &      &          &          & \ddots &          \\[2mm]
1 & d_{0 n}^2  & d_{1 n}^2 & d_{2 n}^2 & \cdots & 0 
\end{array} \right] 
\enspace.
\end{equation} 
The multivariate polynomial $\Gamma_n:=\det(\CM_n)\in\Z[d_{01},d_{02},\ldots,d_{(n-1)\,n}]$
is the \textit{Cayley-Menger determinant}. 

\smallskip 

\par\noindent
Let $v_0, \dots, v_n \in \R^n$ be  a set of $n+1$ points and
denote by $S$ its convex hull in $\R^n$. This determinant gives a formula for the $n$-dimensional volume of $S$ in terms of
the Euclidean distances
$\{{\delta_{ij}}:=\dist(v_i, v_j) \, : \ 0\leq i<j\leq n\}$ 
among these points. 
We have~\cite[~Sec. IV.40]{Blu53}, \cite[~Sec. 9.7]{Ber87}
$$
(-1)^{n+1}\, 2^n \, (n!)^2 \, \Vol_n(S)^2 =  \Gamma_n({\delta_{01}, \delta_{02}, \dots,
\delta_{(n-1)\, n}}) \enspace.
$$
This formula shows that $\Gamma_n$ is a  homogeneous 
polynomial of degree $2\, n$. 
The second polynomial $\Gamma_2$ can be completely 
factorized, giving rise to 
the well-known {\it Heron's formula}
for the area $A$ of a triangle with edge lengths $a$, $b$, and $c$:
\begin{equation} \label{chavo}
16 \, A^2 = 
- \Gamma_2(a,b,c) =  (a+b+c)(-a+b+c)(a-b+c)(a+b-c)
\enspace. 
\end{equation}


\vspace{20mm} 

\begin{figure}[htbp]

\begin{picture}(0,0) 

\put(-60,0){\epsfig{file=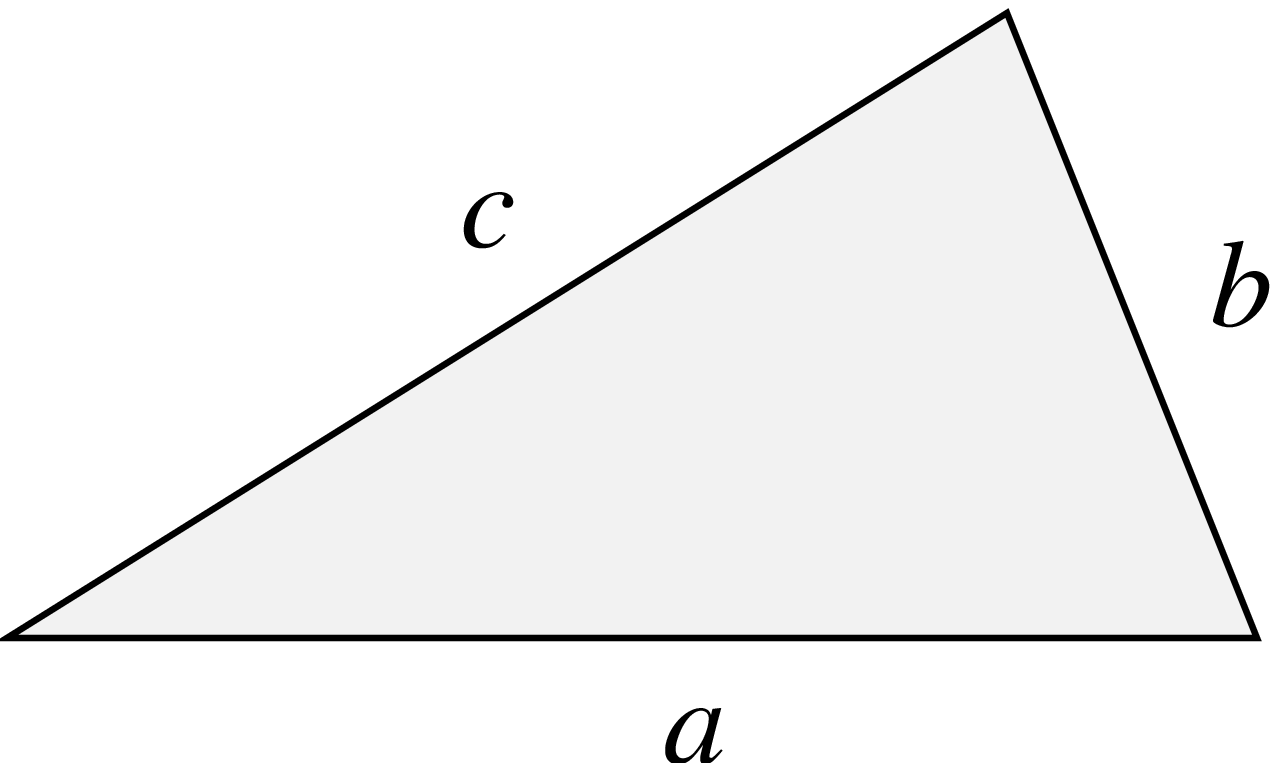, height= 23 mm}} 

\end{picture} 

\end{figure}

\vspace{-5mm} 


Note also that the equation $\Gamma_n(\delta_{01},\delta_{02},\ldots,
\delta_{(n-1)\,n})=0$ gives a necessary and
sufficient condition for the points
$v_0,\ldots,v_n$ to lie in a proper affine subspace of $\R^n.$
\par The Cayley-Menger determinant can be also used for deciding 
whether a set of positive real numbers $\{\delta_{ij}:\,0\leq i<j\leq n\}$
can be realized as the set of edge lengths of an $n$-dimensional simplex in $\R^n$:
in \cite[~Sec. 9.7.3]{Ber87} it is shown that this condition is
equivalent to $(-1)^{n+1}\,\Gamma_n(\delta_{01},\delta_{02},\ldots,
\delta_{(n-1)\,n})>0$.

\smallskip
\par The matrix $\CM_n$ also gives a criterion to determine if $n+2$ points in $\R^n$ lie in an $(n-1)$-dimensional sphere, and         to solve the related problem
of computing the radius of the sphere circumscribed around a simplex. 
To do this, consider the $(1,1)$-minor
$$
\Delta_n:=\det\left[ \begin{array}{ccccc}
0         & d_{01}^2 & d_{02}^2 & \cdots & d_{0 n}^2 \\[2mm] 
d_{01}^2  & 0        & d_{12}^2 & \cdots & d_{1 n}^2 \\[2mm] 
d_{02}^2  & d_{12}^2 & 0        & \cdots & d_{2 n}^2 \\[2mm] 
\vdots &            &          & \ddots &          \\[2mm]
d_{0 n}^2  & d_{1 n}^2 & d_{2 n}^2 & \cdots & 0 
\end{array} \right] \in\Z[d_{01},d_{02},\ldots,d_{(n-1)\,n}]\enspace. 
$$
From this expression we see that this is a homogeneous polynomial 
of degree $2\,n+2$. 

Assume now that $v_0,\ldots,v_n$ do not lie in a proper affine subspace, 
so that $S$ is an $n$-dimensional simplex. The
radius $\rho(S)$ of the sphere circumscribed around $S$ is given by 
\begin{equation}\label{kiko}
\rho(S)^2=-\frac12\frac{\Delta_n(\delta_{01},\delta_{02},\ldots,
\delta_{(n-1)\,n})}{\Gamma_n(\delta_{01},\delta_{02},\ldots,
\delta_{(n-1)\,n})}\enspace.
\end{equation}
Also, the condition for $n+2$ points $v_0,\ldots,v_{n+1}$ in $\R^n$ to lie in 
the same sphere or hyperplane is given by the annulation of the
$(n+1)$-th polynomial
$\Delta_{n+1}(\delta_{01},\delta_{02},\ldots,
\delta_{n\,(n+1)})$ $=0$, see \cite[~Sec. 9.7.3.7]{Ber87}.
\par\smallskip
The third polynomial $\Delta_3$ factorizes  as
\begin{eqnarray} \label{donramon} 
\Delta_3&=&-(d_{01}\, d_{23} + d_{02} \, d_{13} + d_{03} \, d_{12})\, 
(d_{01}\, d_{23} + d_{02} \, d_{13} - d_{03} \, d_{12})\\[0mm] 
&&(d_{01}\, d_{23} - d_{02} \, d_{13} + d_{03} \, d_{12})\,
(-d_{01}\, d_{23} + d_{02} \, d_{13} + d_{03} \, d_{12}) \nonumber 
\enspace.
\end{eqnarray}
This is equivalent to Ptolemy's theorem, 
which states that a convex quadrilateral with edge lengths $a,b,c,d$ 
and diagonals $e,f$ as in the picture, 
is circumscribed in a circle if and only if $a\,c+b\,d=e\,f$.


\vspace{30mm} 

\begin{figure}[htbp]

\begin{picture}(0,0) 

\put(-60,0){\epsfig{file=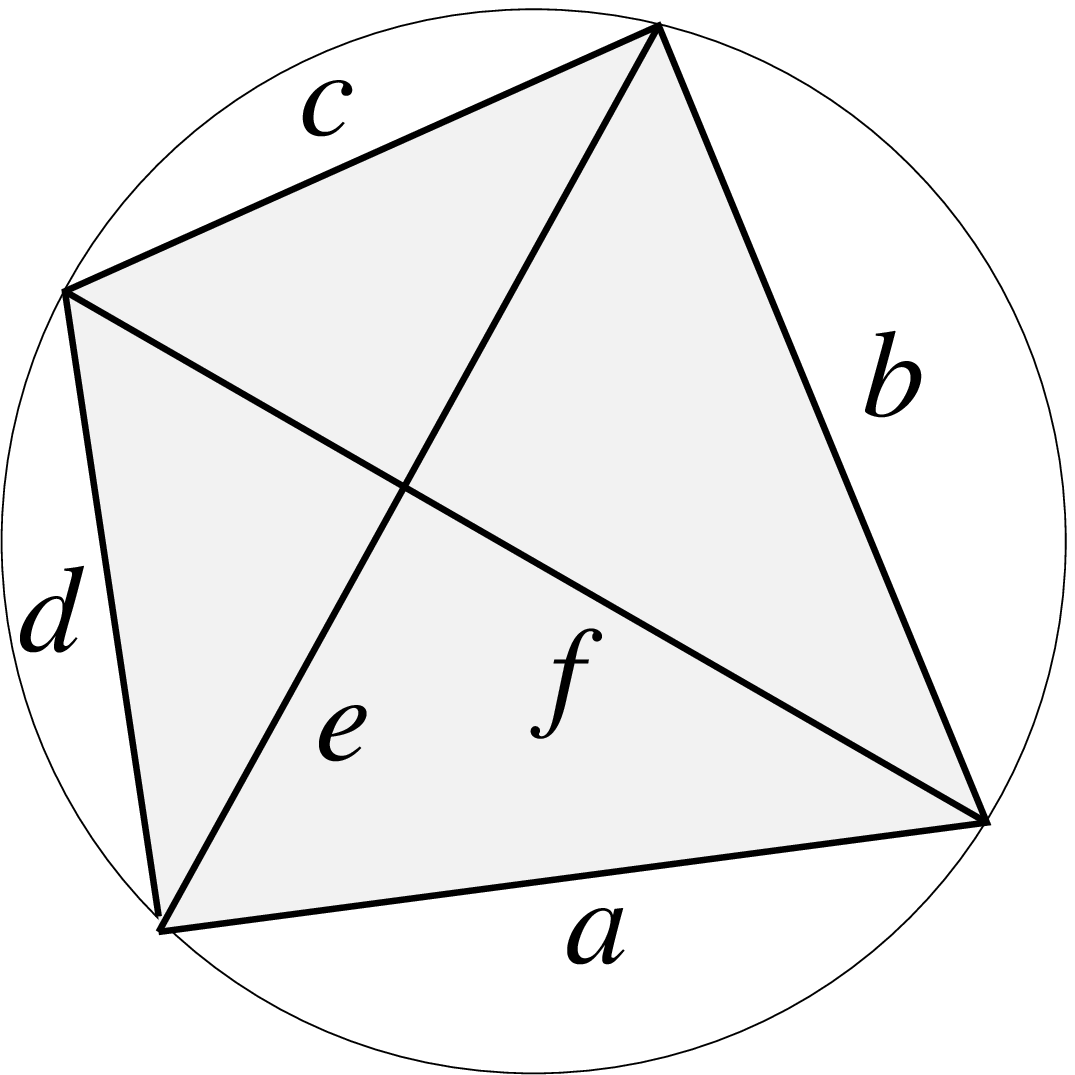, height= 33 mm}} 

\end{picture} 

\end{figure}

\vspace{-3mm} 


The key sources for the Cayley-Menger determinant are the classical 
books by L. Blumenthal \cite{Blu53} and by M. Berger \cite{Ber87}.

\smallskip

This polynomial plays an important role in some problems of metric geometry.
It was first applied by K. Menger in 1928, 
to characterize Euclidean spaces in metric terms alone~\cite[~Ch.~IV]{Blu53}.
It also appears in the metric characterization of 
Riemannian manifolds of constant sectional curvature obtained by M. Berger \cite{Ber81}.

Another important result based on the Cayley-Menger determinant is the proof of the invariance of the volume 
for flexible $3$-dimensional polyhedra
(the ``bellows'' conjecture), see 
\cite{Sab96,CSW97,Sab98}.
There is also a huge literature about 
applications to the study of 
spatial shape of molecules (stereochemistry), see e.g. \cite{KD80,DM00}.

\bigskip 

It is natural to ask whether Heron's formula~(\ref{chavo})
generalizes to higher dimensions, 
that is whether $\Gamma_n$ splits as a product of linear forms. 
Note also that $\Gamma_1= 2 \, d_{01}^2$. 
The purpose of this paper is to prove that this is not possible for $n
\ge 3$. 
Moreover, we show that for $n \ge 3$ 
the only factors 
of $\Gamma_n$ in $\C[d_{0 1}, d_{0 2}, \dots, d_{(n-1)\, n}]$
are the trivial ones, that is either a constant or a constant multiple of
$\Gamma_n$. 
In other words $\Gamma_n$ is {\it absolutely irreducible}.

\begin{thm} \label{irreduc} 
The polynomial \ $\Gamma_n$ 
is irreducible over $\C[d_{0 1}, d_{0 2}, \dots,  d_{(n-1)\, n}]$
for $n\ge~3$. 
\end{thm} 

In a similar way, one may wonder whether $\Delta_n$ splits as a product of
simpler expressions, as in (\ref{donramon}).
Note that $\Delta_1=-d_{01}^4$ and $\Delta_2= 2\, d_{01}^2\,
d_{02}^2\, d_{12}^2$.  
Again we can show that this is not possible for $n\ge 4$. 

\begin{thm} \label{chilindrina} 
The polynomial $\Delta_n$ 
is irreducible over $\C[d_{0 1}, d_{0 2}, \dots,  d_{(n-1)\, n}]$
for $n\ge~4$. 
\end{thm} 

As a straightforward consequence of this, 
we find that the determinant of the general symmetric $n \times n$
matrix  with 
zeros in the diagonal is an absolutely irreducible polynomial for 
$n\ge 4$, see Remark \ref{nionio}. 

\medskip 

We can verify that $\Gamma_3$ is {\it twice} an integral polynomial 
and the same holds for $\Delta_4$. 
This does not affect  their irreducibility over $\C[d_{ij}]$: 
2 is trivial factor as it is a unit of $\C[d_{ij}]$. 
Nevertheless it is interesting to 
determine how they split over $\Z[d_{ij}]$. 
Recall that the {\it content} of an integral polynomial is defined as
the gcd of its coefficients. 

\begin{thm} \label{contenido} 
Let $n \in \N$, then both $\Gamma_n$ and $\Delta_{n+1}$ 
have content 1 for even $n$ and 2 for odd $n$. 
\end{thm}

Let us denote  $\ov{\Z} $ the ring of algebraic integers, that is
the ring formed by  elements in the algebraic closure $\ov{\Q}$ 
satisfying a {\it monic} integral 
equation. 
It is well-known that an integral
polynomial is irreducible over
$\ov{\Z}[d_{ij}]$
if and only if it 
is  irreducible over 
$\C[d_{ij}]$ and has content 1. 
Set 
$$
\begin{array}{ll} 
I_n := \left\{ \begin{array}{ll} 
\Gamma_n&  \hspace{10mm} \mbox{ for {\it n} even} \enspace \\[2mm] 
\Gamma_n/2  & \hspace{10mm}  \mbox{ for {\it  n} odd}  \enspace
\end{array} \right. 
\quad, & \quad \quad
J_n := \left\{ \begin{array}{ll} 
\Delta_n/2&  \hspace{10mm} \mbox{ for  {\it n} even} \enspace \\[2mm] 
\Delta_n  & \hspace{10mm}  \mbox{ for {\it n} odd} \enspace
\end{array} \right. 
\end{array}.
$$
Hence 
Theorems~\ref{irreduc}, \ref{chilindrina} and \ref{contenido} 
can  be equivalently rephrased as the fact that 
$I_n$ and $J_{n}$ are irreducible over 
$\ov{\Z}[d_{ij}]$
(and in particular over ${\Z}[d_{ij}]$) 
for $n\ge 3$ and for $n \ge 4$, respectively.

\bigskip

Let $t_n$ be a new variable and set 
$$
\Lambda_{n,n-1}:= \Gamma_n\big( d_{i n} \mapsto t_n \, : \, 0 \le i \le  n-1 
\big) 
\in \Z[d_{0 1}, d_{0 2}, \dots, d_{(n-2)\, (n-1)}][t_n] \enspace.
$$
Up to a scalar factor, $\sqrt{\Lambda_{n, n-1}}$ is the
formula for the volume of an
isosceles simplex $S(\tau) \subset \R^n$ with base 
$B:= \Conv(v_0, \dots, v_{n-1})$ and 
vertex $v_n$ equidistant at distance $\tau$ to the other 
vertices. 


\vspace{21mm} 

\begin{figure}[htbp]

\begin{picture}(0,0) 

\put(-38,-15){\epsfig{file=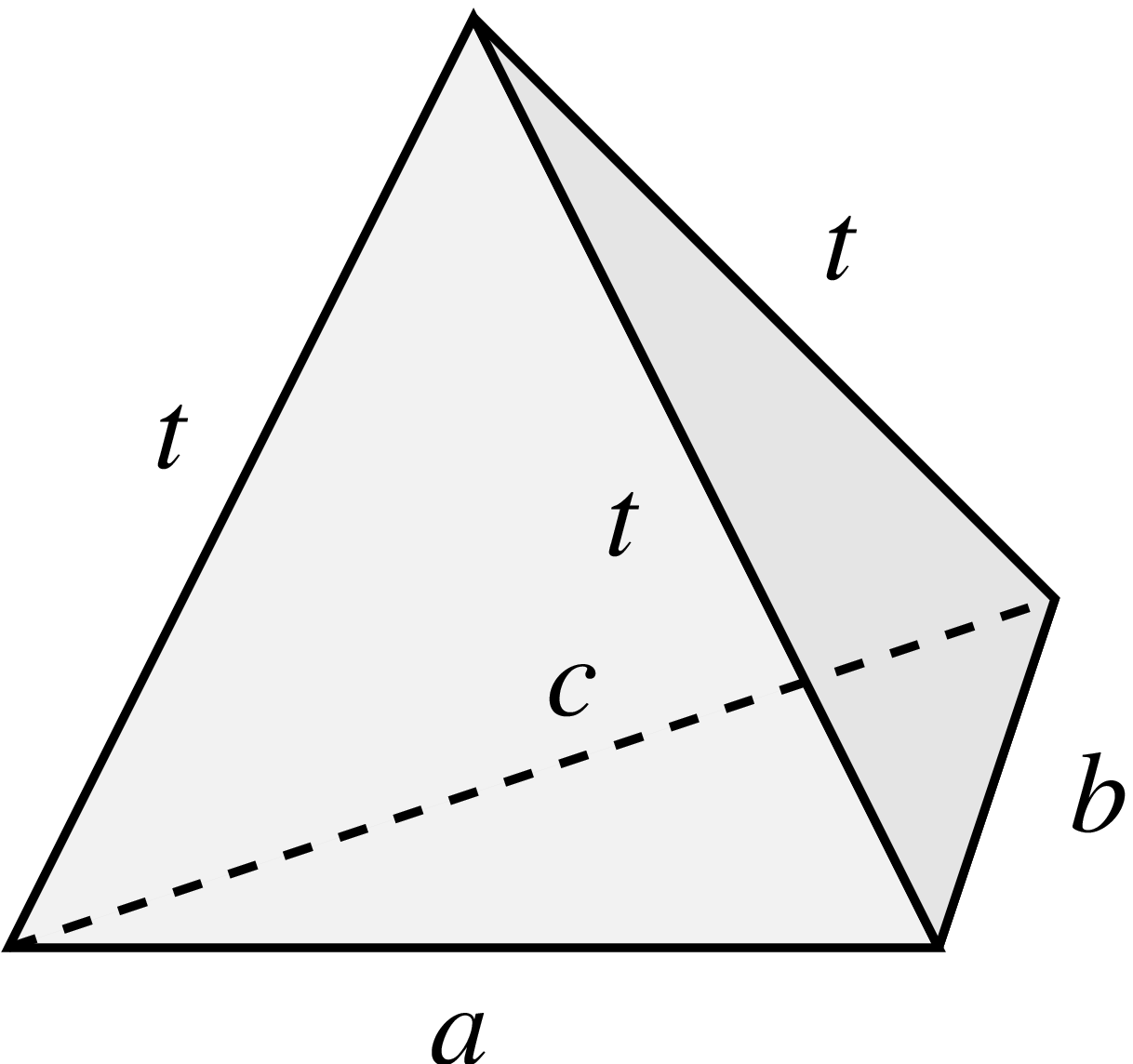, height= 29 mm}} 

\end{picture} 

\end{figure}

\vspace{1mm} 


In~\cite[Sec. 9.7.3.7]{Ber87} it is shown that
\begin{equation} \label{recurrencia} 
\Lambda_{n, n-1}= -2 \, \Gamma_{n-1} \, t_n^2 - \Delta_{n-1} .
\end{equation} 
The dominant term in this expression  corresponds with the geometric 
in\-tui\-tion 
$$ 
\Vol_n(S(\tau)) \sim \frac{\tau}{n} \, \Vol_{n-1} (B) 
\quad \quad \mbox{for }  \tau \gg
0 \enspace.
$$ 
Assuming $\dim(B)=n-1$, note that when 
$\tau=\rho(B)$ is the radius of the circle circumscribing $B$
we have $\Lambda_{n, n-1}=0$ 
and thus we recover (\ref{kiko}).

More generally, let $ 1 \le p \le n$
and set 
$$
\Lambda_{n,p}  := 
\Gamma_n\big( d_{i \ell} \mapsto t_\ell \, : \
p+1 \le \ell \le n, \ 0 \le
i \le \ell-1 \big) 
$$
where $\{t_2 ,\dots, t_n\}$ denotes a further group of variables; 
in particular $\Lambda_{n, n} = \Gamma_n.$ 
It turns out that $\Lambda_{n,p}$ is a homogeneous evaluation of $\Gamma_n$, and
so $\Lambda_{n,p} $ is a homogeneous polynomial 
in $\Z[d_{ij} \, : \ 0 \le i < j \le p ] [t_{p+1}, \dots, t_n]$
of degree $2\,n$. 

\smallskip 

Set $B_p:= \Conv(v_0, \dots, v_p) $ a $p$-dimensional 
simplex with edge
lengths $\{ \delta_{i j} \, : 0 \le i < j \le p \}$, and 
for $0 \ll \tau_{p+1} \ll \cdots \ll \tau_n$ 
set $S(\tau_{p+1}, \dots, \tau_n) \subset \R^n$ the 
$n$-dimensional simplex 
built from $B_p$ by 
successively adjoining a vertex $v_\ell$ equidistant to the previous
vertices $v_0, \dots, v_{\ell-1}$ 
for $\ell = p+1, \dots, n$.
Up to a scalar factor, $ \sqrt{\Lambda_{n,p}}$ is 
the formula for the volume of $S(\tau_{p+1}, \dots, \tau_n)$.  
We have the recursive relation: 

\begin{lem} \label{seniorbarriga} 
$\displaystyle 
\Lambda_{n, p } = -2 \, \Lambda_{ n-1,p}\, t_n^2 - \Lambda_{n-2,p}\,
t_{n-1}^4 $ 
for $ n \ge p+2$. 
\end{lem} 

\begin{proof} 
From the determinantal expression of $\Delta_n$ we get 
\begin{equation} \label{brujadel71}
\Delta_{n-1}(d_{i\, (n-1)}\mapsto t_{n-1} \, : \ 0 \le i \le n-2) = t_{n-1}^4 \,
  \Gamma_{n-2} \enspace, 
\end{equation}  
and so
by (\ref{recurrencia}) we have   $\displaystyle 
\Lambda_{n, n-2 } = -2 \, \Lambda_{ n-1,n-2}\, t_n^2 - \Lambda_{n-2,n-2}\,
t_{n-1}^4 $ 
for $n\ge 2$. 
The general case follows by evaluating $
d_{i \ell} \mapsto t_\ell$ for
$ p+1 \le \ell \le n-2$ and  $  0 \le
i \le \ell-1$ in both sides of this identity. 
\end{proof} 

\medskip 

Theorem \ref{irreduc} is a particular case of the following: 

\begin{prop} \label{doniaflorinda} 
The polynomial \ $\Lambda_{n,p}$ 
is irreducible over $\C[d_{0 1}, d_{0 2}, \dots,  d_{(n-1)\, n}]$
if and only if $n\ge~3$ and $ 2 \le p \le n$. 
\end{prop} 

The following is a graphical visualization of this proposition. 
We encircle the integral points $(n,p)$ such that 
$\Lambda_{n,p}$ is
absolutely irreducible, and we mark with a cross the points where it is not. 
The behavior of $\Gamma_n$ is read from the diagonal.


\vspace{30mm} 

\begin{figure}[htbp]

\begin{picture}(0,0) 

\put(-38,-15){\epsfig{file=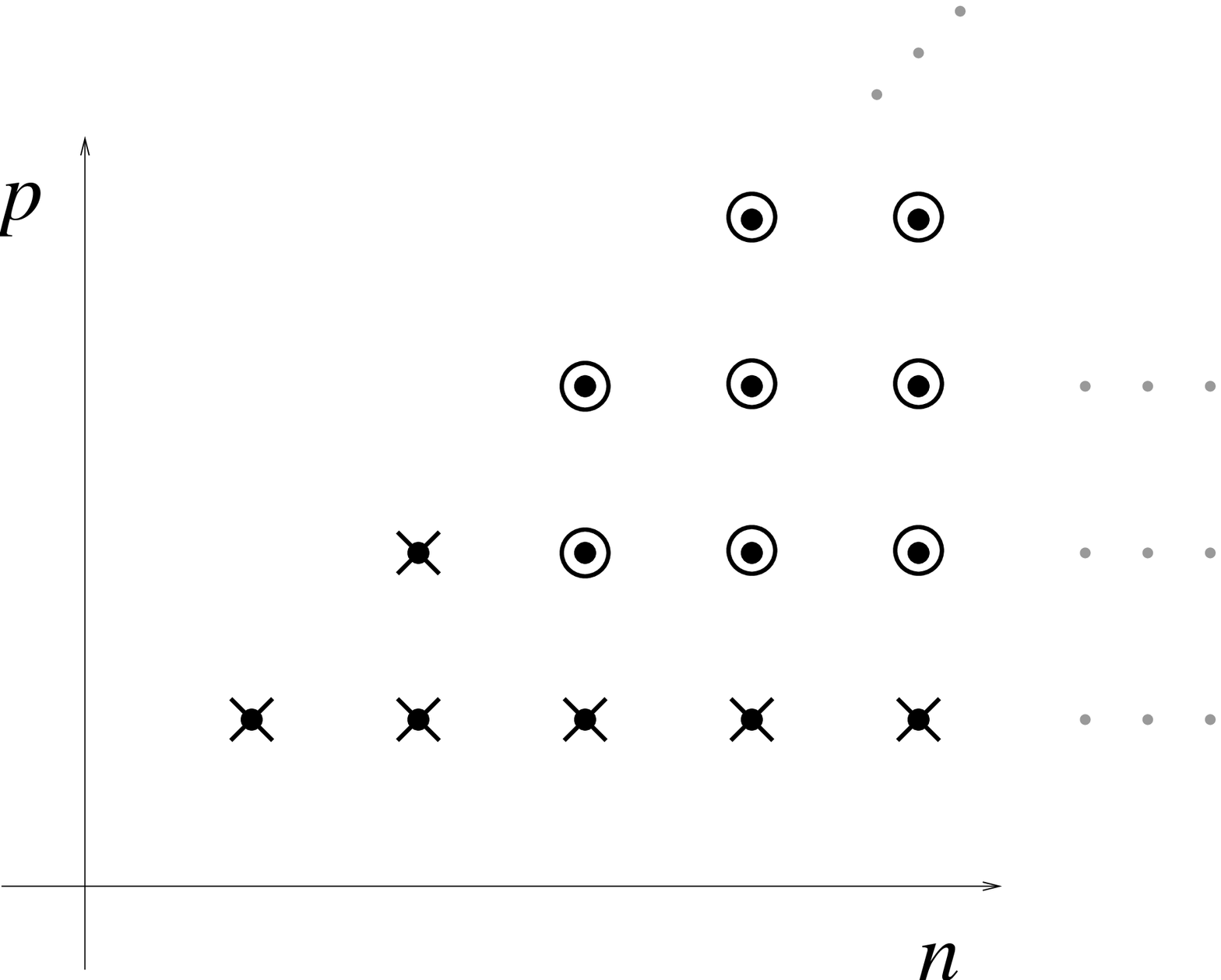, height= 38 mm}} 

\end{picture} 

\end{figure}

\vspace{2mm} 


\begin{proof}
First we will prove by induction that $\Lambda_{n,2}$ is absolutely
irreducible for $n\ge 3$. 
Let $n=3$, and to simplify the notation set 
$d_{01} \mapsto  a, d_{12} \mapsto  b, d_{02} \mapsto c, t_3\mapsto \tau$. 
Identity (\ref{recurrencia}) and Heron's formula imply
\begin{equation} \label{Q3}
\hspace*{10mm} \frac{\Lambda_{3,2}}{2}= - \Gamma_2 \, \tau^2 - \frac{\Delta_2}{2}= 
 (a+b+c)(-a+b+c)(a-b+c)(a+b-c) \, \tau^2 - a^2 \, b^2 \, c^2 \enspace. 
\end{equation} 
The polynomials 
$
f:= (a+b+c)(-a+b+c)(a-b+c)(a+b-c)$ and $g:= - a^2 \, b^2 \, c^2 $
have no common factors over $\C[a,b,c]$ 
and so a (non trivial) factorization of $\Lambda_{3,2}$ should be of the 
form
$$
\Lambda_{3,2} =  (\alpha \, \tau + \beta) ( \gamma \, \tau + \delta)
$$
with $\alpha \,  \gamma = 2\, f$, $\beta \,  \delta = 2\, g$ and $\alpha \, \delta + 
\beta \, \gamma = 0$. 
But this is impossible since  $\alpha,\gamma, \beta\,\delta$ 
have no common factor; we conclude that $\Lambda_{3,2}$ is irreducible. 

\smallskip 

Now let $n \ge 4$ and assume that  $\Lambda_{n-1,2} $ 
is irreducible. 
By Lemma \ref{seniorbarriga} 
$$
\Lambda_{n, 2 } = -2 \, \Lambda_{ n-1,2}\, t_n^2 - \Lambda_{n-2,2}\,
t_{n-1}^4 \enspace.
$$ 
The polynomials $\Lambda_{ n-1,2}$ and $\Lambda_{n-2,2}\,
t_{n-1}^4 $ are coprime, 
since $\Lambda_{n-1,2}$
is irreducible
of degree $2\, n -2$
and  $\deg(\Lambda_{n-2,2}) = 2\, n -4$. 
As before, this implies that any non trivial factorization 
of $\Lambda_{n,2}$ should be of the form
$$
\Lambda_{n,2} = (-2 \, \Lambda_{n-1,2} \, t_n + \beta) ( t_n + {\delta})
$$
with $\beta, \delta \in 
\C[ d_{i   j} \, : 0 \le i < j \le 2 ] [t_{3}, \dots, t_{n-1}] $
such that
$\beta \, \delta = -\Lambda_{n-2,2}\,
t_{n-1}^4 $ and 
$-2\, \Lambda_{ n-1,2} \, \delta + \beta=0$. 
But this is impossible
because $\Lambda_{ n-1,2}$ and  $\Lambda_{n-2,2}\,
t_{n-1}^4$ are coprime. 
We conclude that $\Lambda_{n,2}$ 
is irreducible. 

\smallskip

Now let $n\ge 3$ and $3 \le  p \le n$. 
Suppose that we can write 
$\Lambda_{n,p}= F \cdot G $ with 
$F,G \in  \C[d_{ij} \, : \ 0 \le i <j \le p][t_{p+1} , \dots, t_n]$ 
homogeneous  of degree $\ge 1$.

The evaluation map $d_{k  \ell} \mapsto t$ 
($0 \le k \le \ell, p+1 \le \ell \le n$) 
is  homogeneous and so 
\begin{eqnarray*} 
F'&:=&
F\big(d_{i\ell} \mapsto t_\ell \, : \ p+1 \le \ell \le n, \ 
0 \le i \le \ell-1\big)
\enspace , \\[2mm]  
G'&=& G\big(d_{i\ell} \mapsto t_\ell \, : \ p+1 \le \ell \le n, \  
0 \le i \le \ell-1\big)
\end{eqnarray*} 
are also homogeneous polynomials of degree $\ge 1$, which would 
give a non 
trivial factorization
of $\Lambda_{n,2}$. 
This shows that $\Lambda_{n,p}$ is also irreducible.

\smallskip 

To conclude, we have to verify 
that $d_{01} | \Lambda_{n,1}$ for all $n$, which follows by 
checking that  $\Lambda_{n,1}(d_{01} \mapsto 0) =0$, due to the fact that the 
second and third rows in the matrix defining $\Lambda_{n,1}(d_{01} \mapsto 0)$ coincide. 
The remaining case $n=p=2$ corresponds to Heron's formula. 
\end{proof}

\medskip 

\begin{Proof}{Proof of Theorem \ref{chilindrina}.}
Set 
$$
\Delta_n':= \Delta_n( d_{in}\mapsto 1 \, : \ 1 \le i \le n-1 ) \in 
\Z[d_{01}, d_{02}, \dots, d_{(n-2)\,(n-1)}, d_{0n}]\enspace.
$$
From the determinantal expression of $\Delta_n$ we get 
\begin{equation} \label{lapopis} 
\Delta_n' = 
d_{0n}^4 \, \Gamma_{n-1}\Big(\frac{d_{01}}{d_{0n}}, \dots, 
\frac{d_{0\,(n-1)}}{d_{0n}}, \, d_{12}, d_{13}, \dots, d_{(n-2)\,(n-1)}\Big)
\enspace.
\end{equation} 
Note that the partial degree of $\Gamma_{n-1}$ in the 
group of variables
\begin{equation} \label{profesorgirafales} 
\{ d_{0i} \, : \ 1\le i \le n-1\} 
\end{equation}
is four. 
Hence $\Delta_n'$ is the homogenization
of $\Gamma_{n-1}$ with respect to these variables, 
with $d_{0n}$ as the homogenization variable. 
This follows again from the same determinantal expression.    

Now let $ p,q \in \C[d_{ij}]$ such that 
$
\Delta'_n = F \cdot G$. 
Since $\Delta_n'$ is homogeneous with respect to the variables
(\ref{profesorgirafales}), we have that $F$ and $G$ are also homogeneous with 
respect to this group. 
Now we dehomogenize this identity by setting $d_{0n}\mapsto 1$ and we find 
$$
\Gamma_{n-1} = F( d_{0n}\mapsto 1)\cdot G( d_{0n}\mapsto 1). 
$$
By Theorem~\ref{irreduc}, $\Gamma_{n-1}$ is irreducible for $n \ge 4$, 
which implies that either $F( d_{0n}\mapsto 1) \in \C$ or 
$ G( d_{0n}\mapsto 1) \in \C$. 
This can only hold if $F$ or $G$ is a monomial in $d_{0n}$, but this 
is impossible since  $d_{0n}$ is the homogenization variable. 
We conclude that $\Delta'_n$ is irreducible.

\smallskip 

Now suppose that $\Delta_n$ {\it can} be factorized, and let 
$P, Q \in \C[d_{ij}]$ be homogeneous polynomials of degree $\ge 1$ 
such that $\Delta_n= P\cdot Q$. 
This implies that 
$\Delta_n'= P'\cdot Q'$ with 
$$
P':=P( d_{in}\mapsto 1 \, : \ 1 \le i \le n-1 ) 
\quad , \quad \quad  Q':= Q( d_{in}\mapsto 1 \, : \ 1 \le i \le n-1 ) \enspace. 
$$
Note that $\deg(\Delta_n') = \deg(\Gamma_{n-1}) +4= 2\, n+2$ and so 
$\deg(\Delta_n') = \deg(\Delta_n)$. 
This implies that both $\deg(P')= \deg(P)\ge 1 $ and $\deg(Q')=\deg(Q) \ge 1$, which 
contradicts the irreducibility of $\Delta'_n$. 
Hence $\Delta_n$ is irreducible. 
\end{Proof} 

\medskip

For the proof of Theorem~\ref{contenido} we need an auxiliary result. 
Let $n \in \N$ and $\{ x_{ij} \, : \  1\leq~i<j\leq n\}$ a set of
${(n-1) \, n}/{2}$ variables. 
Then set 
$$
X_n:=\left[
\begin{array}{ccccc}
0&x_{12}&x_{13}&\ldots&x_{1n} \\
x_{12}&0&x_{23}&\ldots&x_{2n} \\
x_{13}&x_{23}&0&\ldots&x_{3n} \\
\vdots&&&\ddots& \\
x_{1n}&x_{2n}&x_{3n}&\ldots&0
\end{array}
\right]
$$
for the general symmetric matrix of order $n$ with zeros in the diagonal.

\begin{lem}\label{mod2}
For odd values of $n$, 
the content of $\det(X_n)$ is divisible by 2. 
\end{lem}

\begin{proof}
Set
$$
A_n:=\left[
\begin{array}{ccccc}
0&x_{12}&x_{13}&\ldots&x_{1n} \\
-x_{12}&0&x_{23}&\ldots&x_{2n} \\
-x_{13}&-x_{23}&0&\ldots&x_{3n} \\
\vdots&&&\ddots& \\
-x_{1n}&-x_{2n}&-x_{3n}&\ldots&0
\end{array}
\right]
$$
for the general {\it antisymmetric} matrix of order $n$. 
Then 
$$
\det(A_n) = \det(A_n^t) = (-1)^n \, \det(A_n) \in \Z[x_{ij}]
\enspace, 
$$
which implies  $\det(A_n) =0$ because $n$ is odd; here $A_n^t$ denotes 
the transpose of $A_n$.
On the other hand 
$ X_n \equiv A_n \pmod{2}$  and so we conclude 
$$
\det(X_n) \equiv \det(A_n) = 0 \pmod{2} \enspace. 
$$
\end{proof}

\begin{Proof}{Proof of Theorem~\ref{contenido}.}
Let $c(n) \in \N$ be the content of $\Gamma_n$. 
Lemma~\ref{mod2} shows that $2|c(n) $ for odd $n$, as the
Cayley-Menger matrix $\CM_n$ 
is symmetric of order $n+2$ with zeros in the
diagonal. 
By Lemma \ref{seniorbarriga} 
$$
\Lambda_{n, n-2} (t_n\mapsto 0) = -   \Gamma_{n-2} \,  t_{n-1}^4 \enspace.  
$$
By definition  
$ \Lambda_{n, n-2}(t_n\mapsto 0) $ is an integral
evaluation of $\Gamma_n$ and so $c(n)$ divides its content, that is 
$c(n) | c(n-2)$. 
We conclude by induction, checking  the statement directly for $n=1$ and
$n=2$.

\smallskip 

Now let $c'(n) \in \N$ be the content of $\Delta_n$. 
Lemma~\ref{mod2} shows that $2|c'(n) $ for even $n$, as the
matrix in the definition of $\Delta_n$ 
is symmetric of order $n+1$ with zeros in the
diagonal. 
Identity (\ref{lapopis}) implies that $c'(n) |c(n-1)$, that is 
$c'(n)=1$ for $n$ odd and $c'(n)|2$ for $n$ even; we conclude that 
$c'(n)=2$ in this case.  
\end{Proof}

\begin{rem} \label{nionio} 
Set
$$
K_n := \left\{ \begin{array}{ll} 
\det(X_n) &  \hspace{10mm} \mbox{ for {\it n} even} \enspace, \\[2mm] 
\det(X_n)/2  & \hspace{10mm}  \mbox{ for {\it  n} odd}  \enspace.
\end{array} \right.  
$$
As a byproduct of Theorems~\ref{chilindrina} and~\ref{contenido}, we find that 
$K_n$ is an irreducible polynomial over $\overline{\Z}[x_{ij}\, : \ 1\le i <j\le n]$
for $n\ge 5$;  
a direct verification shows that this is also true for $n=4$.

\end{rem}


\subsection*{Acknowledgments} 
It is a pleasure to thank Victor Alexandrov from the Sobolev Institute of Mathematics
at Novosibirsk, who posed us this problem and constantly encouraged us to solve it.
We also thank him for pointers to the litterature on the Cayley-Menger
determinant. 
 
We are grateful to Mat{\'\i} as Gra{\~n}a for suggesting us a 
simplified proof of Lemma~\ref{mod2}.

\smallskip

Part of this paper was written while Mart{\'\i} n Sombra was visiting the University of California at Berkeley, supported by the Miller Institute for Basic
Research in Science. 


\typeout{Referencias}


\begin{thebibliography}{XXXXX}

\addcontentsline{toc}{section}{References}

{\small

\bibitem[Ber81]{Ber81}
Berger, Marcel,
\newblock {\em  Une caract{\'e}risation purement m{\'e}trique des vari{\'e}t{\'e}s riemanniennes {\`a} courbure constante.\/} 
\newblock E. B. Christoffel (Aachen/Monschau, 1979),  pp. 480--492, Birkh{\"a}user, Basel-Boston, Mass., 1981.
 
\bibitem[Ber87]{Ber87}
Berger, Marcel, 
\newblock{\em Geometry I.\/} 
\newblock Springer-Verlag, Berlin, 1987.

\bibitem[Blu53]{Blu53}
Blumenthal, Leonard M., 
\newblock{\em Theory and applications of distance geometry.\/} 
\newblock Oxford, at the Clarendon Press, 1953.

\bibitem[CSW97]{CSW97}
Connelly, R.; Sabitov, I.; Walz, A., 
\newblock{\em The bellows conjecture.\/}  
\newblock Beitr{\"a}ge Algebra Geom.~{\bf 38}  (1997) 1--10. 

\bibitem[DM00]{DM00}
Deo, Narsingh; Micikevicius, Paulius,
\newblock{\em Generating edge-disjoint sets of quadruples in parallel for the molecular conformation problem.\/} 
\newblock Congr. Numer.  {\bf 143}  (2000) 81--96.

\bibitem[KD80]{KD80}
Klapper, Michael H.; DeBrota, David,
\newblock {\em Use of Cayley-Menger determinants in the calculation of molecular structures.\/} 
\newblock J. Comput. Phys.  {\bf 37}  (1980) 56--69. 

\bibitem[Sab96]{Sab96}
Sabitov, I. Kh., 
\newblock{\em The volume of a polyhedron as a function of its metric. 
(Russian).\/}  
\newblock Fundam. Prikl. Mat.  {\bf 2}  (1996) 1235--1246.

\bibitem[Sab98]{Sab98}
Sabitov, I. Kh.,
\newblock{\em The volume as a metric invariant of polyhedra.\/}
\newblock Discrete Comput. Geom. {\bf 20} (1998) 405--425.

}

\end{thebibliography}
\end{document}